\documentclass[11pt]{amsart}
\usepackage{latexsym,amssymb,amsmath,amsthm}
\usepackage{comment}
\usepackage{float}
\usepackage{graphicx}
\newtheorem{thm}{Theorem}

\theoremstyle{definition}

\newcommand{\Z}{\mathbb{Z}}

\def\G{\mathcal{G}}
\def\M{\mathcal{M}}

\def\tr{\rm{tr}}
\def\str{\rm{str}}

\title{An integrable evolution equation in geometry}
\author{Oliver Knill}
\date{April 30, 2013}
\address{
        Department of Mathematics \\
        Harvard University \\
        Cambridge, MA, 02138
        }

\begin{document}
\maketitle
\begin{abstract}
We introduce an integrable Hamiltonian system which Lax deforms the Dirac operator 
$D=d+d^*$ on a finite simple graph or compact Riemannian manifold.
We show that the nonlinear isospectral deformation always leads to an expansion of the original space,
featuring a fast inflationary start. The nonlinear evolution leaves the Laplacian $L=D^2$ invariant
so that linear Schr\"odinger or wave dynamics is not affected. The expansion has the following effects: 
a complex structure can develop and the nonlinear quantum mechanics asymptotically becomes 
the linear relativistic Dirac wave equation $u''=Lu$. While the later is not aware of the expansion
of space and does not see the emerged complex structure, nor the larger non-commutative geometric setup, the 
nonlinear flow is affected by it. The natural Noether symmetries of quantum mechanics introduced here 
force to consider space as part of a larger complex geometry. 
The nonlinear evolution equation is a symmetry of quantum mechanics which still features supersymmetry, 
but it becomes clear why it is invisible: while the McKean-Singer formulas $\str(e^{i D(t) t}) = \str(e^{-L t})=\chi(G)$ 
still hold, the superpartners $f,Df$ are orthogonal only at $t=0$ and become parallel or anti-parallel for 
$|t| \to \infty$. 
\end{abstract}

\section{A Lax pair} 

The Dirac operator $D=d+d^*$ on a compact Riemannian manifold or on a finite simple graph \cite{knillmckeansinger} 
is the square root of the Hodge Laplacian $L=d d^* + d^* d$, where $d$ is the exterior derivative. Together 
with $B=d-d^*$ satisfying $B^2=-L,\{B,D \; \}=0$ we can look at the  Lax pair \cite{Lax1968}
\begin{equation}
D'=[B,D] \; . 
\label{lax}
\end{equation}
This is motivated by \cite{Kni94} and of course \cite{Toda}. While the Laplacian has nonnegative spectrum, 
the Dirac operator $D$ has the symmetric spectrum $\pm \lambda_j$. 
The deformed operator has the form $D(t) = d+d^*+b$. In the simplest case,
the differential equations $d' = d b-b d, b' = d d^* - d^* d$
preserve $d^2=(d^*)^2=0$ and $L=\{d,d^* \; \}$. We also have $\{d,b \; \}=\{d^*,b \; \}=0$. 
If $df=0$ is a cocycle, and $f'=b(t)f(t)$ then $d(t) f(t)=0$ so that $f(t)$ remains a
cocycle. If $f=dg$ is a coboundary and $g'=bg$ then $f(t) = d(t) g(t)$ so that $f(t)$ remains a coboundary. 
The system therefore deforms cohomology: graph or de Rham cohomology does not change if we use $d(t)$ instead
of $d$. Define the unitary operator $U(t)$ by $U'=B U, U(0)=1$. From $\{B,D \;\}=0$ we get $D'=2BD$. 
As for any Lax pair, since $d/dt ( U^* D U)=0$, we see that $D(t) = U D U^*$ is isospectral to $D(0)$. 
Because $L'=D'D+D D'=(BD-DB)D + D(BD-DB)=[B,L]=[B,-B^2]=0$,
the Laplacian $L$ is time independent. The solution $u(t)=\cos(Dt) u(0) + \sin(Dt) D^{-1} u'(0)$ of
the linear wave equation $u''=Lu$ only involves the Laplacian $L$ and therefore is not affected by the deformation.
The nonlinear wave $u(t) = \cos(D(t) t) u(0) + \sin(D(t) t) D^{-1}(0) u'(0) = {\rm Re}( e^{i D(t) t} \psi(0) )$
with $\psi(0) = u(0) - i D(0)^{-1} u'(0)$ however depends on $D$. This time-dependent Hamiltonian flow will
in the complex be replaced by the much more natural integrable nonlinear flow $U(t) u(0)$ satisfying $u'(t) = B(t) u(t)$ where
$B(t)$ is like $i D(t)$ a square root of $-L$.  
Since $d(t)$ is a deformed exterior derivative, we can look at the new Dirac operator
$C(t) = d(t) + d(t)^*$ and the new Laplacian $M(t) = C(t)^2$ and $V(t) = b(t)^2$. 
We check that the Laplacian $L=D^2=M+V=C^2+b^2$ 
decomposes into two commuting operators. 
We have $B'=[D,b]$.  The operators $R=d d^*$ and $S=d^* d$ commute because 
their products $RS,SR$ are zero: initially, $[b,d d^*]=0$, we differentiate to get $[b',d d^*]=0$. We also
have $[b,(d d^*)' = [b,(db-bd) d^* + d (b d^* -d^*b)] = [b,-bd d^*-d d^* b]=0$
showing that the four operators $b,R=d d^*,S=d^* d,V=b^2$ commute and can be simultaneously diagonalized. 
With $P(x)=(-1)^{p(x)}$, with $p(x)=p$ on $p$ forms $\Omega^p$, 
we have the supersymmetric relations 
$\{P,C \; \}=\{P,B \;\}=0,\{P,P \; \}=1,\{C,C \; \}=M,\{B,B \;\}=-M$ which describe the individual
deformed geometries. As for any exterior derivative $d_i: \Omega^i \to \Omega^{i+1}$, the anti-commutation relations
$\{d_i,d_j\} = \delta_{ij} L_i$ hold, where $L_p: \Omega^p \to \Omega^p$ is the Laplacian on $p$-forms. 
Supersymmetry fails for the deformed operator $D(t)$ in the sense that  $\{D,P \; \} \neq 0$ and $D$ does no more produce an
isomorphism between the Bosonic subspace $\oplus \Omega^{2k}$ and Fermionic subspace $\oplus \Omega^{2k+1}$ of the
exterior bundle $\Omega = \oplus \Omega^k$. The symmetry is still present however because $P(t)$ could be deformed also
but we can not save the invariant splitting into Fermions and Bosons. 
A key observation for the scattering analysis is to see that 
$O=b d d^*$ is selfadjoint with no negative eigenvalues and $Q=b d^* d$ is selfadjoint
with no positive eigenvalues. Proof:
the eigenvalues of $O,Q$ are initially zero because $b$ is zero at $t=0$. Look at an eigenvalue $\lambda$ 
of $O$. From the Rayley formula $\lambda' = \langle w O',v \rangle$,
where $v$ is the unit eigenvector and $w$ the dual vector $(O^*)' w = \lambda w$ and by symmetry, we have 
$\lambda' = \langle  v,O'v \rangle$. We know that $b$ and $d d^*$ have the same eigenvectors for $t>0$
because they commute. Using $(d d^*)' = -2 b d d^*$ we get
$O' = (d d^*)^2 + b (d d^*)' = (d d^*)^2 - 2 b^2 d d^* = (d d^*)^2 - 2 b O$. 
If $v$ is an eigenvector to $O$ to the eigenvalue $\lambda$, then 
$\lambda' = \langle v,O'v \rangle = \langle v,((d d^*)^2 v -2b O) v\rangle$. This means that 
$\lambda' \geq 0$ if $\lambda=0$. In other words, eigenvalues can not cross $0$. 
The computation for $Q$ is similar. We have Lyapunov functions because $\tr(b^2)' \geq 0$ and so $\tr(M') \leq 0$:
$d/dt \tr(b^2) = 2 \tr(b b')=2 \tr(b d d^* - b d^* d) \geq 0$. 
Initially, we have $\tr(b^2)'=0$ and asymptotically, we have $\tr(b^2)=\tr(L)$
so that $\tr(b^2)'$ will have a maximum somewhere. We see an initial inflation both in the positive and
negative time direction. We have $\tr(bC) = \tr(Cb)=0$ because these matrices do not have anything in the diagonal.
From $\tr(L)'=0$ follows $\tr(M(t))' =\tr(L-C b - b C - b^2)'= -\tr(b^2)' \leq 0$. 
We see that $M(t)$ has its spectrum $\sigma(M(t)) \in [0,a(t)]$ with $a(t) \to 0$ for $|t| \to \infty$.
It follows that $d(t)$ and $C(t)$ converge to zero and $b(t)$ converges to an operator 
$V$ satisfying $V^2=L$. If $\tr(M')=0$ then $b d d^*=0$ which is not possible for $t>0$ by looking at eigenplanes.
For every $k>0$, we have $\str(B^k)=0$ at $t=0$. This is clear for odd $k$, because there is
nothing in the diagonal. For even $k$, note that $B^2=L$ at all times, so that $B^{2n}=L^n$. But we still know
that $\str(L^n)=0$, by the classical McKean-Singer formula \cite{McKeanSinger,knillmckeansinger}.
The nonlinear analogue $\str(U(t)) = \chi(G)$ of McKean-Singer holds:
we know that $\str(L^n)=0$ initially because of the linear McKean Singer result. Because
$\str(U(t))$ is real analytic and $\str(U(0))=\chi(G)$ by definition,
it is enough to verify that $d^k/dt^k \str(U) = 0$ for all $k$ at $t=0$. To see this,
differentiate $U$ at $t=0$ to get $U'=BU, U''=B' U + B^2 U=B U^2+B^2 U, U'''=B U^3+B^2 U^2 + 2 B^2 U + B^3$.
Using $U(0)=Id$ and $\tr(B^k)=0$ one can deduce that all these derivatives are zero.
We have $D(t) + D(-t) = 2C(t)$ because this is true at $t=0$ and because $b(t)+b(-t)=0$ for all $t$. 
In particular, the attractors satisfy $D(\infty) = -D(-\infty)$. We have sketched the proof of:

\begin{thm}
Both in the graph and manifold case, system (\ref{lax}) has the property that the limit 
$\lim_{t \to \infty} D(\pm t) = \pm D(\infty)$ and $\lim_{t \to \infty} U(t)$ exists. 
We have $D(t) = U^*(t) D(0) U(t)$ for some unitary $U$ satisfying the McKean-Singer equations 
$\str(U(t))=\chi(G)$ for all $t$.  
The Laplacian $L=D^2$ does not change but it is for all $t$ a sum $L=M+V$ of two commuting operators
where $M=C^2$ is a new Laplacian for a new Dirac operator $C(t)=d(t)+d(t)^* \to 0$ belonging to a new exterior 
derivative which has the same cohomology than $d$, and a block diagonal part $V=b^2$ whose limit is a square root of $L$. 
\end{thm}

\section{Examples}

{\bf 1)} In the Riemannian manifold case, $D$ is a differential operator on the exterior bundle of $M$. The deformed
operator is a pseudo differential operator. The deformation can be described by matrices,
once a basis is chosen. For the circle $M=T$, the Dirac operator is
$D = \left[ \begin{array}{cc} 0           & \partial_x \\ -\partial_x  & 0 \end{array} \right]$
and $L = D^2= \left[ \begin{array}{cc} -\Delta  &  0 \\ 0      &  -\Delta \end{array} \right]$
leaves 0-forms and 1-forms invariant. Fourier theory gives an eigenbasis 
$\left[ \begin{array}{c} \pm i e^{i n x} \\ e^{i n x} \end{array} \right]$
belonging to eigenvalues $\pm n$ so that $\sigma(D)$ is the set of integers $\Z$ and $L$ 
has the eigenvalues $n^2$ for $n=0,1,\dots$. The zeta function of the circle 
$\zeta(s) = \sum_{n \neq 0} n^{-s} = \zeta_{MP}(s/2) + (-1)^s \zeta_{MP}(s/2) = \zeta_{MP}(s/2) [ 1+e^{i \pi s}]$ is 
analytic in the entire complex plane as for any odd dimensional Riemannian manifold. 
The factor $1+e^{i \pi s}$ has naturally regularized the poles of 
the Minakshisundaram-Pleijel zeta function $\zeta_{MS}$. For the circle, the Dirac zeta function
has the same roots then the classical Riemann zeta function. 
Since $\zeta'(0) = -1$, the circle has the Dirac Ray-Singer determinant ${\rm det}(D)=e$. 
The deformation $D=D(t)= d+d^* + b = \left[ \begin{array}{cc} B           & A                    \\
                                                              A^*         & C         \end{array} \right]$
satisfies $D' = [B,D]$ using 
$B = d-d^*= \left[ \begin{array}{cc} 0                     & A                    \\
                                    -A^*                   & 0          \end{array} \right]$ can be written as the matrix 
differential equation
\begin{equation}
B' = 2 A A^*, A' = 2 A C,  C' = -2 A^* A \; . 
\label{2}
\end{equation}
Since the quantity $AC+BA$ is time invariant,
$L$ is block diagonal with entries $B^2 + A A^*$ and $C^2 + A^* A$. 
In Fourier space, $A,B,C$ are double infinite matrices. We have $B(0)=C(0)=A(\infty)=0$ and
$A(0) = {\rm Diag}(\dots -3i,-2i,-i,0,i,2i,3i, \dots)$ and 
$B(\infty) = -C(\infty) = {\rm Diag}( \dots 3,2,1,0,1,2,3, \dots )$. 
System (\ref{2}) shows initial inflation and asymptotic exponential expansion of the individual circles. 
Also for general initial conditions $A,B,C \in M(n,C)$, system~(\ref{2}) 
satisfies $A(t) \to 0$. 
The deformed Dirac operator $D(t)$ describes a larger geometry of two separate circles.  \\
{\bf 2)} For a finite simple graph $G$ with simplex set $\G_k$ of cardinality $v_k=|\G_k|$, Euler characteristic
$\chi(G) = \sum_{k=0}^{\infty} (-1)^k v_k$, the Dirac operator $D$ is a $v \times v$ matrix with 
$v=\sum_{k=0}^{\infty} v_k$. For the two point graph $G=K_2$, with $v_0=2,v_1=1$ the Dirac operator is
the $3 \times 3$ matrix $D=\left[ \begin{array}{ccc} 0 & 0 & -1 \\ 0 & 0 & 1 \\ -1 & 1 & 0 \end{array} \right]$
with eigenvalues $-\sqrt{2},\sqrt{2},0$. The Laplacian is
$L= D^2= \left[ \begin{array}{ccc} 1 & -1 & 0 \\ -1 & 1 & 0 \\ 0 & 0 & 2 \\ \end{array} \right]$.
The zeta function is $\zeta(s) = 2^{-s/2} + (-2)^{-s/2}$, the pseudo determinant $2$. 
Starting with $D(0) =  \left[ \begin{array}{ccc}  0 &  0 & -1 \\ 0 &  0 &  1 \\
                                   -1 &  1 &  0 \end{array} \right]$ we end up with 
$D(\infty) = -D(-\infty) =   \left[ \begin{array}{ccc}
                      1         &   -1      &  0   \\
                     -1         &    1      &  0   \\
                      0         &  0        &  -2   \end{array} \right]/\sqrt(2)$. 
The differential equation for
$D = \left[ \begin{array}{ccc}
                      b         & c         & -d        \\
                      c         & b         &  d        \\
                     -d         & d         &  e       \end{array} \right] \; ,  
   B = \left[ \begin{array}{ccc}
                      0         & 0         & -d        \\
                      0         & 0         &  d        \\
                      d         & -d        &  0       \end{array} \right]$
simplify to
\begin{equation}
b'=2d^2, d'=-4bd 
\label{3}
\end{equation}
because $c=-b$ and $e=-2b$. Explicit solutions of (\ref{3}) are
$d(t)=(1-\tanh^2(\sqrt{8} \; t))^{1/2}$, $b(t)=\tanh(\sqrt{8} \; t)/\sqrt{2})$
and the integral $d^2+2b^2$ whose level curves are ellipses. Looking at
$d'(t)=-4 d(t) b(t)$ confirms initial inflation for the increasing distance of 
the two points $x,y$.

\section{Remarks}
1) {\bf Hamiltonian formalism.} 
To see this system as a Hamiltonian system to the Hamiltonian $\tr(L^2)$ in the manifold case, we have to 
consider Dirac zeta function $\zeta(s)=\sum_{\lambda \neq 0} \lambda^{-s}$. Note however
that to define and run the flow, we do not need the Hamiltonian.  The 
traces $\zeta(-k)=\tr(D^k)$ and the determinant $-\zeta'(0)$ are then defined by analytic continuation. 
Unlike the zeta function of the Laplacian \cite{LapidusFrankenhuijsen} which is only meromorphic, we can chose the branch 
$\zeta(s)=(1+e^{i \pi s}) \zeta_L(s/2)$, where $\zeta_L$ is the zeta function of the Laplacian, so that the Dirac zeta function 
has an analytic extension everywhere both in the graph or compact Riemannian manifold case. 
The branch is natural because it 
leads to the trivial roots $\tr(D^{2k+1}=\zeta(-2k-1)$ agreeing with the fact that the operator $D^{2k+1}$ has only $0$'s in the diagonal. 
For every observable $H=\tr(L^n)$ and for even positive $n$, there is a Hamiltonian 
flow $D'=[B,h'(D)] = J \nabla H(D)$.  But because $h'(D)=D p(L)$ and because $L$ commutes with $B$, higher degree flows 
are related by the first flow with Hamiltonian $\tr(D^2)$ by a energy dependent 
time change on each eigenplane.  \\

2) {\bf Broken supersymmetry.} The eigenfunctions $f,Df$ are initially perpendicular
superpartners and become parallel or anti-parallel for $|t| \to \infty$, hiding supersymmetry:
if $f$ is a Boson, then $D(t) f$ is only a Fermion at $t=0$.
$C(t)$ defines a Connes pseudo metric $\sup_{|Df|_{\infty}=1} |f(x)-f(y)|$ \cite{Connes} leading to an expanding 
space which features a fast inflationary start near the supersymmetric origin. $D$ determines via
non-commutative geometry a pseudo distance in the larger space
$\G = \bigcup G_p$, the union of the set $\G_p$ of $p$-cliques or
the exterior bundle $\M = \bigcup M_p$ which is a union of p-branes $M_p$, images of the 
Gelfand transform of the Banach algebra of $p$-forms with pointwise multiplication. 
While the linear Dirac evolution preserves supersymmetry, the nonlinear evolution provides
a mechanism to see it broken for measurements. \\

3) {\bf Distances}. $G_k$ always has infinite distance to $G_0$ 
in the Connes pseudo metric $d$ because $|Df|=0$ constant $f \in \Omega_0$. More generally, if the cohomology group
$H^p(G)$ is nontrivial then $\G_p$ has infinite distance to all other $\G_q$ for $p \neq q$ and furthermore,
$d$ is not a metric. An other distance in $\G$ or $\M$ is obtained by measuring
how long a wave $\psi=u+iv$ with $|v|=1$ and $u$ localized at $x$ needs to get to $y$. 
if the minimal $t$ such that $e^{i D t} (x+iv) = y+iw$ for $|v|=1$, then $d(x,y)=1/t$. 
Also in the nonlinear case, define $d(x,y)$ to be the minimal $t$ such that for some $|v|=1$ we have $U(t)(x+iv) = (y+iw)$. 
This is possible if ${\rm Re}( U(t) (x+iv)) = y$ has a solution $|v|=1$. Linear or nonlinear quantum mechanics makes it
possible that with $\G$ or $\M$, different k-branes $\G_k$ or $\M_k$ have finite distance. \\

4) {\bf Deformed curvature.} Gauss-Bonnet-Chern \cite{cherngaussbonnet} and Poincar\'e-Hopf \cite{poincarehopf} 
and its link \cite{indexexpectation,indexformula} are the key to define curvature on the larger
space. For $t=0$, the spaces $\G_p$ or $\M_p$ are disconnected and curvature $k(x)$ is the expectation 
$K(x)={\rm E}[i_f(x)]$ of the index $i_f(x)$ defined for functions on $\G_p$ or $\M_p$. At $t=0$ we take
the product measure $P_0$. The unitary evolution $U(t)$
pushes forward this probability measure $P_0$ on functions. The new measure $P_t$ produces a new expectation
and so deforms curvature $K_t(x)$.
We have not yet explored the question how curvature defined by the new operator $C(t) = d(t) + d(t)^*$ of the expanding space 
evolves if we rescale length so that the diameter stays constant. 
The answer will be interesting in any case: if the evolution simplifies space, it can be used as a 
geometric tool which unlike the Ricci flow features global existence.
If geometry is not simplified, it should lead to limiting geometries, or a dynamical system on an 
attractor of limiting geometries. \\

\begin{figure}[H]
\parbox{15.4cm}{
\parbox{5cm}{\scalebox{0.08}{\includegraphics{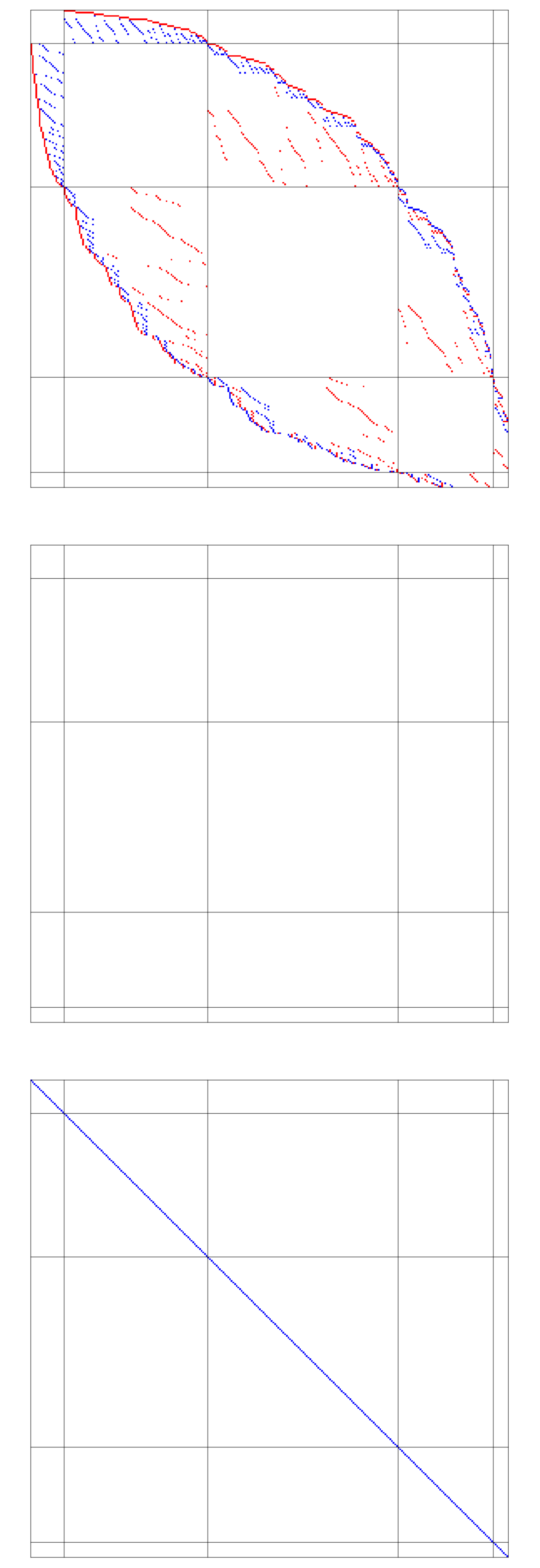}}}
\parbox{5cm}{\scalebox{0.08}{\includegraphics{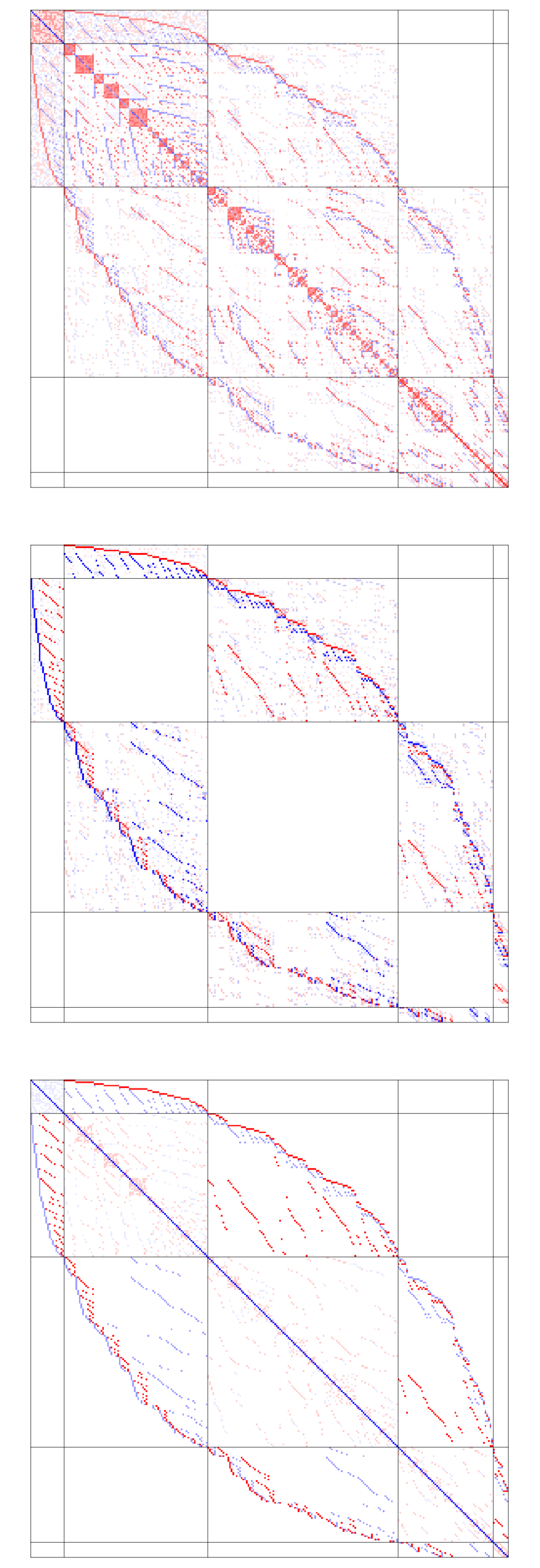}}}
\parbox{5cm}{\scalebox{0.08}{\includegraphics{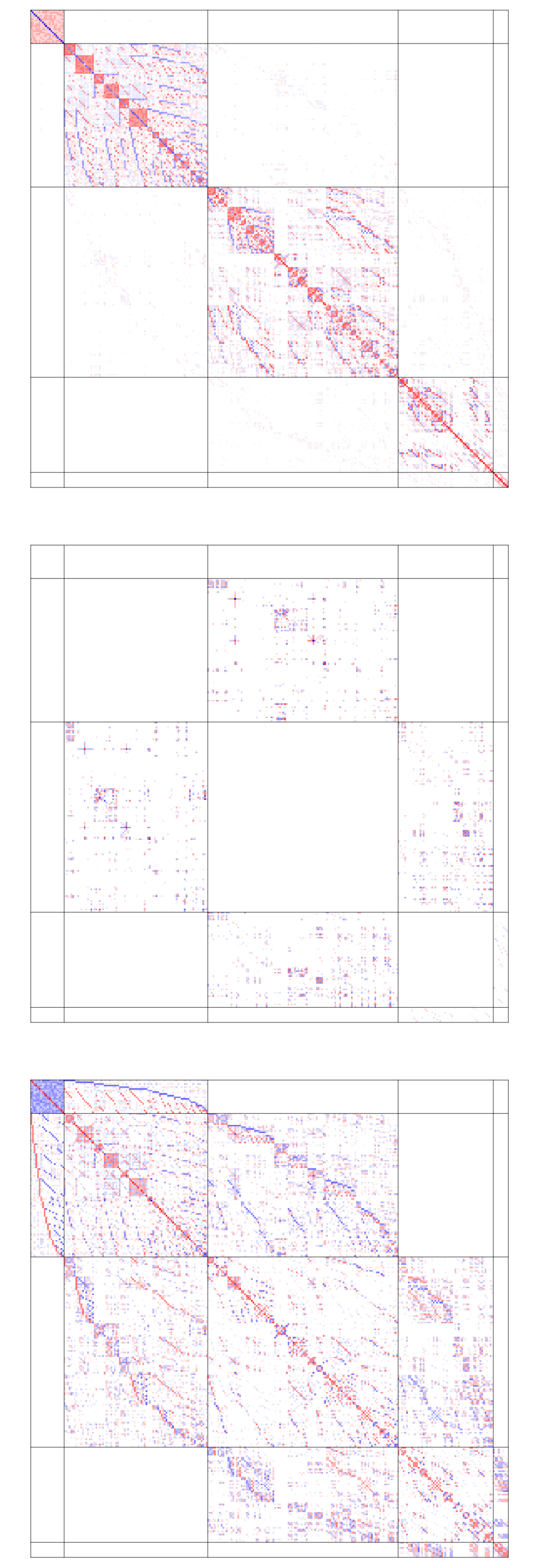}}}
}
\caption{
The complex evolution of a Dirac operator $D$ of a randomly chosen graph $G$ with 20 vertices, 86 
edges, the Poincar\'e polynomial encoding the Betti numbers is $p(t)=1 + 3t + 2t^2$, 
the clique polynomial is $c(t)=20 + 86t + 114t^2 + 57t^3 + 9t^4$ and
$\chi(G)=p(-1)=c(-1)$ by Euler-Poincar\'e.
We see snapshots of the $v(1) \times v(1)$ matrices ${\rm Re}(D(t)),{\rm Im}(D(t)),{\rm Re}(U(t))$ 
at $t=0,t=0.2,t=1$. The complex part appears and disappears. 
\label{figure1}
}
\end{figure}

5) {\bf Complex structure.} A complex generalization of the system is obtained by defining 
$B=d-d^*+i \beta b$, where $\beta$ is a parameter. This is similar to \cite{Toda} who modified the Toda flow
by adding $i \beta$ to $B$. The case $\beta=0$ is now a very special case. 
For $\beta=1$, we have $B^2=-C^2-b^2=-L$ and $A=iD$ satisfies $A^2=-L$ but no more $\{B,D \;\}=0$ as 
in the real case. The essential features of the system like expansion are $\beta$-independent because $b(t)$
and so $M(t)$ is $\beta$ independent. The limit $D(\infty)$ for the complex flow is the same than the 
limiting real flow. The Lax pair $D'=[B,D]$ is now more symmetric and the 
nonlinear unitary flow $U(t)$ is asymptotic to the linear Dirac flow because $B(t) \to b$, where $b$ 
is a square root of $L$. Remarkably, a complex 
differential structure has emerged during the evolution because $D(t)$ is complex for $t>0$ even so we have start 
with a real graphs or manifold. Define $\partial={\rm Re}(d)$ and 
$\overline{\partial}=i{\rm Im}(d))/2$, then $\partial^2 = \overline{\partial}^2 = 0$ and 
$d=\partial + \overline{\partial}$. Because cocycles and coboundaries deform in an explicit way, 
cohomology groups defined by $\partial$ and $\overline{\partial}$ are both the same than for $d$.
Since $\partial \overline{\partial} = \overline{\partial} \partial=0$, the Laplacian $L=D^2$ is
the sum of two Laplacians $L^{\partial} = (D^{\partial})^2$ and 
$L^{\overline{\partial}} = (D^{\overline{\partial}})^2$, where $D^{\partial} = \partial + \partial^*$
and $D^{\overline{\partial}} = \overline{\partial} + \overline{\partial}^*$. One could call a graph 
$G$ with a complex structure given by $D$ a K\"ahler graph,
if $D^{\partial} = D^{\overline{\partial}}$. 
The complex structure disappears asymptotically for $t \to \infty$ 
(see Figure~(\ref{figure1})) justifying
that  $\exp(i D(t) t)$ is close to $U(t)$ satisfying $U(t)^* D(t) U(t)=D(0)$. 
While the exterior derivative $d$ as well as the new Dirac operator $C=d+d^*$ are complex for $t>0$, the
operator $M=C^2$ is real if we start with a real $D$. While asymptotically, 
$||{\rm Im}(D(t))||/||{\rm Re}(D(t))|| \to 0$, the complex structure is especially relevant in the 
early stage of the evolution. 

\begin{figure}[H]
\parbox{16.9cm}{
\parbox{9.5cm}{\scalebox{0.25}{\includegraphics{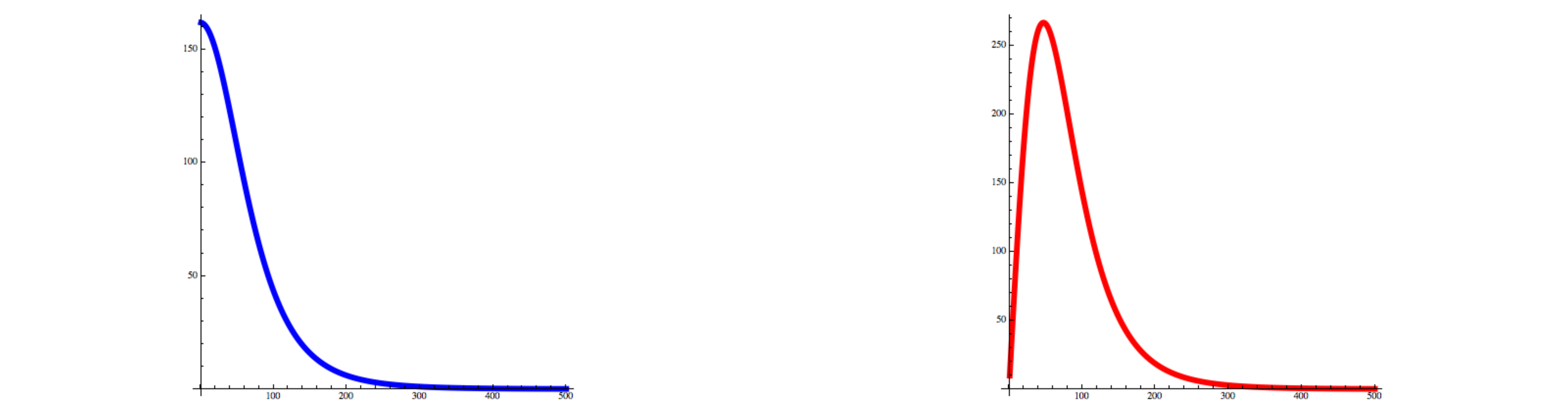}}}
}
\caption{
\label{figure2}
The evolution of the function ${\rm tr}(M(t))$ and $\frac{d}{dt} {\rm tr}(M(t))$
showing inflation. This figure is again done using a random Erdoes-Renyi graph $G$. 
It looks the same for all graphs. 
}
\end{figure}


\end{document}